\documentclass[12pt,english,reqno]{amsart}
\usepackage{amssymb,amsmath,amstext,amsthm,amsfonts}
\usepackage[dvips]{graphicx}
\usepackage[ansinew]{inputenc}
\usepackage{amssymb}

\usepackage[english,francais]{babel}

\newtheorem{theorem}{Theorem}[section]
\newtheorem{lemma}[theorem]{Lemma}

\newtheorem{corollary}[theorem]{Corollary}

\def\og{\leavevmode\raise.3ex\hbox{$\scriptscriptstyle\langle\!\langle$~}}
\def\fg{\leavevmode\raise.3ex\hbox{~$\!\scriptscriptstyle\,\rangle\!\rangle$}}

\def \RR {{\mathbb R}}
\def \ZZ {{\mathbb Z}}
\def \NN {{\mathbb N}}

    \def \co {\mathcal{O}}

\newcommand{\leb}{\operatorname{Leb}}

\newcommand{\qand}{\quad\mbox{and}\quad}

\begin{document}


\selectlanguage{english}


\vspace{-2.6cm}

\selectlanguage{francais}
\title{Backward volume contraction  for \newline \   endomorphisms with eventual volume
expansion}


\selectlanguage{english}
\author{J. F. Alves}
\address{Centro de Matem\'atica da Universidade  do Porto\\
Rua do Campo Alegre 687, 4169-007 Porto, Portugal}
\email{jfalves@fc.up.pt}
\author{Armando Castro}
\address{Departamento de Matem\'atica, Universidade Federal da Bahia\\
Av. Ademar de Barros s/n, 40170-110 Salvador, Brazil}
\email{armando@impa.br}
\author{Vilton Pinheiro}
\address{Departamento de Matem\'atica, Universidade Federal da Bahia\\
Av. Ademar de Barros s/n, 40170-110 Salvador, Brazil}
\email{viltonj@ufba.br}

\thanks{Work carried out at the  Federal University of
Bahia. Partially supported by FCT through CMUP and UFBA}

\maketitle

\begin{abstract}
We consider smooth maps on compact Riemannian manifolds. We prove
that under some mild condition of eventual volume expansion
Lebesgue almost everywhere  we have uniform backward volume
contraction  on every pre-orbit of Lebesgue almost every point.

\vskip 0.5\baselineskip

\selectlanguage{francais}  \noindent{\bf R\'esum\'e} \vskip
0.5\baselineskip \noindent Nous consid\'erons des transformations
diff\'erentiables sur des variet\'es Riemannienes compactes. Nous
montrons que dans une certaine condition mod\'er\'ee d'expansion
de volume nous pouvons d\'eduire  que pour Lebesgue presque chaque
point nous avons contraction uniforme de volume en arri\`ere  de
chaque pr\'e-orbite.
\end{abstract}

\selectlanguage{english}

\section{Statement of results}

Let $M$ be a compact Riemannian manifold and let $\leb$ be a
volume form on $M$ that we call Lebesgue measure. We take $f\colon
M\to M$ any smooth map.
 Let $0<a_1\le a_2 \le a_3\le \dots$ be a sequence  converging
to infinity. We define
 \begin{equation}\label{aga}
    h(x)=\min\{n>0 \colon |\det Df^n(x)|\ge a_n\},
\end{equation}
if this minimum exists, and $h(x)=\infty$, otherwise. For $n\ge
1$, we take
 \begin{equation}\label{gamma}
   \Gamma_n=\{x\in M \colon h(x) \ge n\}.
\end{equation}

\begin{theorem} \label{theo} Assume that $h\in L^p(\leb)$, for some $p>3$, and
take $\gamma<(p-3)/(p-1)$. Choose any sequence $0<b_1\le b_2 \le
b_3\le \dots$ such that $b_kb_n\ge b_{k+n}$ for every $k,n\in
\NN$, and assume that there is $n_0\in\NN$ such that $b_n\le
\min\left\{a_n,\leb(\Gamma_n)^{-\gamma}\right\} $ for every $n\ge
n_0$. Then, for $\leb$ almost every $x\in M$, there exists $C_x>0$
such that $| \det Df^n(y)|>C_x b_n$ for every $y\in f^{-n}(x).$
 \end{theorem}

\smallskip

We say that $f\colon M\to M $ is {\em eventually volume expanding}
if there exists $\lambda>0$ such that for Lebesgue almost every
$x\in M$
 \begin{equation}
 \sup_{n\ge 1}\frac1n\log|\det Df^n(x)|> \lambda.
 \end{equation}
Let $h$ and $\Gamma_n$ be defined as in~(\ref{aga})
and~(\ref{gamma}), associated to the sequence $a_n=e^{\lambda n}$.

\medskip

 \begin{corollary} \label{Mtheo} If $f$ is eventually volume expanding, then for Lebesgue almost
every $x\in M$ there are $C_x>0$ and $\sigma_n\to \infty$ such
that $| \det Df^n(y)|>C_x \sigma_n $ for every $y\in f^{-n}(x)$.
Moreover, given $\alpha>0$ there is $\beta>0$ such that
\begin{enumerate}
\item   if $\leb(\Gamma_n)\le \co(e^{-\alpha n})$, then we may
take $ \sigma_n \ge e^{\beta n}$;
    \item if
$\leb(\Gamma_n)\le \co(e^{-\alpha n^\tau})$ for some $\tau>0$,
then we may take $\sigma_n \ge e^{\beta n^\tau}$;
    \item if
$\leb(\Gamma_n)\le \co(n^{-\alpha})$ and $\alpha>2$, then we may
take $\sigma_n \ge n^\beta$.
\end{enumerate}
 \end{corollary}

 \medskip

Specific rates  will be obtained in Section~\ref{se.examples} for
some eventually volume expanding endomorphisms. In particular,
non-uniformly expanding maps such as quadratic maps and Viana maps
will be considered.

For the proof of our results we give abstract versions of the
techniques developed by Armando Castro in his PhD. thesis \cite{C0}
and articles (\cite{C1}, \cite{C2}). More precisely, we adapt his
chain concatenation ideas and Redundance Elimination Algorithm to
noninvertible contexts. This is main target in the next section.

%

\section{Concatenated collections}

 Let  $(U_n)_n$ be a collection of measurable subsets
 of $M$ whose union covers a full Lebesgue measure subset of  $M$. We say that
$(U_n)_n$ is a {\em concatenated collection}
  if:
  $$x\in U_n \qand f^n(x)\in U_m\quad\Rightarrow\quad x\in
  U_{n+m}.$$
Given $x\in \bigcup_{n\ge 1} U_n$, we define
 $u(x)$ as the minimum
$n\in\NN$ for which $x\in U_n$. Note that by definition we have
$x\in U_{u(x)}$. 
We define  the {\em chain generated by $x\in \bigcup_{n\ge 1}
U_n$} as
 $C(x)=\{x,f(x),\dots,f^{u(x)-1}(x)\}.$

\smallskip

Just as Corollary 2.9 in \cite{C2}, the next Lemma is a consequence
of Borel-Cantelli Lemma. It says that, if we see the collection of
chains given by a concatenated collection $(U_n)_n$ as a tower, and
this tower has finite measure, than a.e. point in $M$ is contained
in just a finite number of chains.

\begin{lemma}
\label{lemma1} Let $(U_n)_n$ be a
 concatenated collection. If
$$
\sum_{n\ge 1}\sum_{j=0}^{n-1}\leb(f^j(u^{-1}(n)))<\infty,
$$ then
we have $ \sup\left\{\,u(y)\ \colon \; y\in \bigcup_{n\ge 1}
U_n\;\mbox{and}\,\; x\in C(y)\,\right\}<\infty $ for Lebesgue
almost every $x\in M$.
\end{lemma}

\smallskip

 Assume that for a given $x\in M$ there exists an infinite
number of chains
 $C_j=\left\{y_j,f(y_j), \dots ,f^{s_j-1}(y_j)\right\}$, $ j\ge 1$, containing $x$
 with $s_j\to\infty$. For each $j\ge1$ let
$1\le r_j<s_j$ be such that $x=f^{r_j}(y_j)$.
 First we verify that
$\lim r_j=\infty$. If not, then replacing by a subsequence, we may
assume that there is $N>0$ such that $r_j<N$ for every $j\ge1$. This
implies that $y_j\in\bigcup_{i=1}^{N}f^{-i}(x)$ for every $j\ge1$.
Since $\#(\bigcup_{i=1}^{N}f^{-i}(x))<\infty$ and the number of
chains is infinite, we have a contradiction. Since $r_j\to\infty$
and $x=f^{r_j}(y_j)\in f^{r_j}(u^{-1}(s_j))$, we have
$x\in\bigcup_{n\ge k}\bigcup_{j=0}^{n-1}f^j(u^{-1}(n))$ for every
$k\ge 1$. Since we are assuming $\sum_{n\ge
1}\sum_{j=0}^{n-1}\leb(f^j(u^{-1}(n)))<\infty$, we have
$\leb\big(\bigcup_{n\ge k}\bigcup_{j=0}^{n-1}f^j(u^{-1}(n))\big)\to
0, $ when $k\to\infty$. This completes the proof of
Lemma~\ref{lemma1}.

\smallskip

If a point $x$ is contained just in a finite number of chains, as in
the Lemma above, it is obvious that the greatest length of such
chains is bounded by $N= N(x)$. Such points have a interesting
property.

\begin{lemma}
\label{lemma2} Let $(U_n)_n$ be a concatenated collection.

If $\sup\left\{\,u(y)\ \colon \; y\in \cup_{n\ge 1}
U_n\;\mbox{and}\,\; x\in C(y)\,\right\}\le N$, 
then $f^{-n}(x)\subset U_{n}\cup \dots \cup U_{n+N}$ for all
$n\ge1$.
\end{lemma}

\smallskip

The idea of the proof  is inspired in Prop. 2.15-2.19 in \cite{C1}
and Prop. 2.12-2.14 in \cite{C2}. Due to the concatenation property,
one can glue chains from $z= f^{-n}(x)$ up to some moment less then
$n$, when for the first time we glue a chain containing $x$. Such
chain can not have length greater then $N$, and therefore $z \in U_n
\cup \dots U_{n+N}$. Let us write down these arguments in detail.

 Assume that $\sup\left\{\,u(y)\ \colon \; y\in \cup_{n\ge 1}
U_n\;\mbox{and}\,\; x\in C(y)\,\right\}\le N$ and take $z\in
f^{-n}(x)$.  Let $z_j=f^j(z)$ for each $j\ge 0$. We distinguish
the cases $x\in C(z)$ and $x\notin C(z)$. If $x\in C(z)$, then
$n\le u(z)\le n+N$. Hence
 $z\in U_{u(z)}\subset U_n\cup\cdots\cup U_{n+N}.$
If $x\notin C(z)$, then  letting $u_0=u(z)$ we must have $u_0<n$.
Let $u_1=u(z_{u_0})$.  If $u_0+u_1< n$ we take
$u_2=u(z_{u_0+u_1})$. We proceed in this way until we find the
first $s\le n $ such that $n\le u_0+ \dots +u_s$. Note that
$u_s=u(z_{u_0+\cdots +u_{s-1}})$, and by the choice of $s$ we must
have $x\in C(z_{u_0+\cdots +u_{s-1}})$. Our assumption implies
that $u(z_{u_0+\cdots +u_{s-1}})\le N$, and so $u_0+ \dots +u_s\le
n+N$. By construction we have
$$
        z   \in   U_{u_0},\,\,
     f^{u_0}(z)=z_{u_0}  \in  U_{u_1} ,
 f^{u_0+u_1}(z)=z_{u_0+u_1}  \in  U_{u_2},\,
 \dots,
 $$
 $$
 f^{u_0+\cdots u_{s-1}}(z)=z_{u_0+\cdots u_{s-1}}\in U_{u_s}
 $$
%
By the definition of a  concatenated collection we conclude that
$z\in U_{u_0+u_1+\dots+u_s}$.

\section{Proofs of main results}\label{se.general}

Let us now prove Theorem \ref{Mtheo}. Suppose that $h\in
L^p(\leb)$, for some $p>3$. This implies that
$\sum_{n\ge1}n^p\leb(h^{-1}(n))<\infty$, and so there exists some
constant $K>0$ such that $\leb(h^{-1}(n))\le Kn^{-p}$ for every
$n\ge1$. Now, taking  $0<\gamma<(p-3)/(p-1)$ we have for some
$K'>0$
$$\sum_{n=1}^{\infty}n \left(\sum_{k=
n}^{\infty}\leb(h^{-1}(k))\right)^{1-\gamma}\le
\sum_{n=1}^{\infty}n (K'/n^{p-1})^{1-\gamma}
<\infty.$$
 Defining $U_n=\{x\in M\ \colon |\det Df^n(x)|\ge
b_n\},$ then we have that $\{ U_1, U_2, \dots \}$ is a
concatenated collection with respect to the Lebesgue measure.
Moreover, setting  $U^*_n= U_n\setminus(
 U_1\cup...\cup U_{n-1})$ one observes that
$U^*_n\subset \cup_{m\ge n}h^{-1}(m)$, for otherwise there would
be $x\in U^*_n\cap h^{-1}(m)$ with $m<n$, and so $a_m\ge b_m>|\det
Df^m(x)|\ge a_m,$ which is not possible. As $|\det Df^j(x)|< b_j$
for every $x\in U^*_n$ and $j<n$, we get $\leb(f^j(U^*_n))\le b_j
\leb(U^*_n)$ for each $j<n$. Hence
\begin{eqnarray*}
\sum_{n=n_0+1}^{\infty}\sum_{j=0}^{n-1} \leb(f^j(U^*_n))&\le&
\sum_{n=n_0+1}^{\infty}\sum_{j=0}^{n-1}b_j \leb(U^*_n)\\
&\le&
\sum_{n=n_0+1}^{\infty}\sum_{j=0}^{n_0-1}b_j \leb(U^*_n)+\sum_{n=n_0+1}^{\infty}\sum_{j=n_0}^{n-1}b_j \leb(U^*_n)\\
&\le&
\sum_{j=0}^{n_0-1}b_j+\sum_{n=n_0+1}^{\infty}\sum_{j=n_0}^{n-1}b_j
\leb(U^*_n)
\end{eqnarray*}
Now, we just have to check that the last term in the sum above is
finite. Indeed,
\begin{eqnarray*}
\sum_{n=n_0+1}^{\infty}\sum_{j=n_0}^{n-1}b_j \leb(U^*_n) &\le&
\sum_{n=n_0+1}^{\infty}\sum_{j=n_0}^{n-1}b_j\sum_{k=
n}^{\infty} \leb(h^{-1}(k))\\
&\le&\sum_{n=n_0+1}^{\infty}n b_n\sum_{k= n}^{\infty} \leb(h^{-1}(k))\\
&\le&\sum_{n=n_0+1}^{\infty}n \left(\sum_{k=n}^{\infty}
\leb(h^{-1}(k)\right)^{-\gamma}\sum_{k=
n}^{\infty} \leb(h^{-1}(k))\\
&=&\sum_{n=n_0+1}^{\infty}n \left(\sum_{k= n}^{\infty}
\leb(h^{-1}(k))\right)^{1-\gamma}<\infty.
\end{eqnarray*}
%
Applying Lemmas~\ref{lemma1}~and~\ref{lemma2}, we get for each
generic point $x\in M$ a positive integer number $N_x$ such that if
$y\in f^{-n}(x)$ then $y\in U_{n+s}$ for some $0\le s\le N_x$.
Therefore, $|\det Df^{n+s}(y)|>b_{n+s}\ge b_{n}$. Taking
$C_x=K^{-N_x}$, where $K=\sup\{|\det Df(z)|\colon  z\in M\},$ we
obtain Theorem~\ref{theo}:
$$|\det Df^{n}(y)|=\frac{|\det Df^{n+s}(y)|}{|\det
Df^{s}(x)|}>C_x b_{n}.$$ 

Now we explain how we use Theorem~\ref{theo} to prove
Corollary~\ref{Mtheo}. Recall that in Corollary~\ref{Mtheo} we have
$a_n=e^{\lambda n}$ for each $n\in \NN$.
 Assume first that $\leb(\Gamma_n)\le \co(e^{-c^\prime n})$ for
some $c^\prime>0$. Then it is possible to choose $c>0$ such that
 $b_n =e^{cn},$ for $n\ge n_0$.
 The other two  cases
are obtained under similar considerations.


\section{Examples: non-uniformly expanding maps}\label{se.examples}

An important class of dynamical systems where we can immediately
apply our results is the class of non-uniformly expanding
dynamical maps introduced in \cite{ABV}. As particular examples of
this kind of systems we present below one-dimensional quadratic
maps and the higher dimensional Viana maps.

\paragraph*{Quadratic maps.}

Let $f_a\colon [-1,1]\to [-1,1]$ be given by $f_a(x)=1-ax^2$, for
$0<a\le 2$. Results in \cite{BC1,J} give that for  a positive
Lebesgue measure set of parameters $f_a$ in non-uniformly
expanding. Ongoing work \cite{F} gives that for a positive
Lebesgue measure set of parameters there are $C,c>0$ such that
$\leb(\Gamma_n)\le Ce^{-c n}$ for every $n\ge1$.

Thus, it follows from  Corollary~\ref{Mtheo} that {\em there
exists $\beta>0$ such for Lebesgue almost every $x\in I$ there is
$C_x>0$ such that $| (f^n)'(y)|>C_x e^{\beta n}$ for every $y\in
f^{-n}(x)$.}
\paragraph*{Viana maps.}

  Let $a_0\in(1,2)$ be such that the critical point $x=0$
is pre-periodic for the quadratic map $Q(x)=a_0-x^2$. Let
$S^1=\RR/\ZZ$ and $b:S^1\rightarrow \RR$ given by $b(s)=\sin(2\pi
s)$. For fixed small $\alpha>0$, consider the map $\hat f$ from
$S^1\times\RR$ into itself given by $\hat f(s, x) = \big(\hat
g(s),\hat q(s,x)\big)$,
 where  $\hat q(s,x)=a(s)-x^2$ with
$a(s)=a_0+\alpha b(s)$, and $\hat g$ is the uniformly expanding
map of $S^1$ defined by $\hat{g}(s)=ds$ (mod $\ZZ$) for some
integer $d\ge2$. 
For $\alpha>0$ small enough there is an interval $I\subset (-2,2)$
for which $\hat f(S^1\times I)$ is contained in the interior of
$S^1\times I$. Thus, any map $f$ sufficiently close to $\hat f$ in
the $C^0$ topology has $S^1\times I$ as a forward invariant
region. Moreover, there are $C,c>0$ such that $\leb(\Gamma_n)\le
Ce^{-c\sqrt n}$ for every $n\ge1$; see \cite{AA,BST,V}.

Thus, it follows from  Corollary~\ref{Mtheo} that {\em there
exists $\beta>0$ such for Lebesgue almost every $X\in S^1\times I$
there is $C_X>0$ such that $| \det Df^n(Y)|>C_X e^{\beta\sqrt n}$
for every $Y\in f^{-n}(X)$.}


\end{document}